\documentstyle[11pt,amstex,amscd]{article}
\def\bold{\bf}
\def\al{\alpha}
\def\b{\betta}

\def\text#1{{\em #1}}

\def\la{\lambda}

\def\ep{\varepsilon}

\def\be{\begin{equation}}
\def\ee{\end{equation}}
\def\bear{\begin{eqnarray}}
\def\eear{\end{eqnarray}}
\def\best{\begin{eqnarray*}}
\def\eest{\end{eqnarray*}}

\newtheorem{theorem}{Theorem}[section]

\newtheorem{lemma}[theorem]{Lemma}
\newtheorem{cor}[theorem]{Corollary}

\newtheorem{defn}[theorem]{Definition}

\def\rem{ \addtocounter{theorem}{1}
{\non \bf Remark \arabic{section}.\arabic{theorem} }}

\def\non{\noindent}

\def\ra{\rightarrow}

\def\r#1{\right#1}
\def\l#1{\left#1}
\def\ma#1{\mathop {#1} \limits}

\def\b{\beta}

\def\ti{\times}

\def\Z{{ \Bbb Z}}

\def\P{{ \Bbb P}}
\def\Q{{ \Bbb Q}}
\def\cx{{ \Bbb C}}
\def\ev{\mbox{\rm ev}}

\def\wt#1{\widetilde{#1}}

\def\ov#1{\overline{#1}}

\def\w{\omega}

\def\M{{\cal M}}

\title{\bf Gromov-Witten Invariants of  Symplectic  Sums\vskip.2in}
\author{ Eleny-Nicoleta Ionel\thanks{both authors partially supported by
the N.S.F.} \\  M.I.T.\\
 Cambridge, MA  02139 \and Thomas H. Parker{\footnotesize *}\\ Michigan
State University\\ East Lansing,
MI   48824}

\date{\today}

\addtocounter{section}{0}
\begin{document}

\maketitle

\vskip.15in

The natural  sum operation for symplectic manifolds is defined
by gluing  along
codimension two submanifolds.  Specifically, let $X$ be a  symplectic
$2n$-manifold with a symplectic $(2n-2)$-submanifold  $V$.  Given a similar
pair $(Y,\overline{V})$ with a
symplectic identification $V=\overline{V}$ and a complex anti-linear
isomorphism between the normal bundles of $V$
and $\overline{V}$, we can form the symplectic sum
$Z=X\#_{V=\overline{V}}Y$.  This note announces a general formula
for computing the Gromov-Witten invariants of the  sum $Z$ in terms of relative Gromov-Witten
invariants of
$(X,V)$ and $(Y,\overline{V})$.

Section 1 is a review of the GW invariants for symplectic manifolds and the
associated
invariants, which we call $TW$ invariants, that count reducible curves.
The corresponding relative
invariants of a symplectic pair $(X,V)$ are defined in section 2.  The  sum
formula is stated in a
special case in section 3, and in general as Theorem 4.1.  The last section
presents two applications: a
short derivation of the Caporaso-Harris formula [CH], and new proof that
the rational enumerative
invariants of the rational elliptic surface are given by the ``modular
form'' (\ref{BL}).

 Related results, involving symplectic sums along contact manifolds,  are
being developed by Li and Ruan
[LR]  and by Eliashberg and Hofer.

\vskip.3in

%%%%%%%%%%%%%%%%%%%%%%%%%%%%%%%%%%%%%%%%%%%%%%%%%%%%%%%%%%%%
%%%%%%%%%%%%%%%%%  Section 1 %%%%%%%%%%%%%%%%%%%%%%%%%%%%%%
%%%%%%%%%%%%%%%%%%%%%%%%%%%%%%%%%%%%%%%%%%%%%%%%%%%%%%%%%%%%%%

\setcounter{equation}{0}
\section{Symplectic Invariants}
\medskip

The  moduli space of $(J,\nu)$-holomorphic
maps from  genus $g$ curves  with $n$ marked points representing a class
$A$ in the free part of
$H_2(X;\Z)$  has a compactification $\overline{\M}_{g,n}(X,A)$.  This
comes with a map
\begin{equation}
\overline{\M}_{g,n}(X,A)\ \longrightarrow \ \overline{\M}_{g,n}\ti
X^n\ti H_2^{\mbox{\scriptsize {free}}}(X;\Z),
\label{modulimap1}
\end{equation}
where the first factor is the ``stabilization'' map $st$ to the
Deligne-Mumford moduli space (defined by collapsing all unstable components
of the domain curve), the second factor records the images of the marked
points, and the last one
keeps  track of the homology class $A$ of the image.  After
perturbation (cf. [LT]), the image defines a rational homology class,
\begin{equation}
\left[\overline{\M}_{g,n}(X,A)\right] \ \in\
H_*(\overline{\M}_{g,n};\Q)\otimes H_*(X^n;\Q)\otimes
H_*(H_2^{\mbox{\scriptsize {free}}}(X;\Z);\Q)
\label{modulimap2}
\end{equation}
for each choice $g,n, A$. The last term in (\ref{modulimap2}) can be
identified with the rational group ring
$RH_2(X)$ of $H_2^{\mbox{\scriptsize{free}}}(X;\Z)$, that is, with finite
sums $\sum c_A
t_A$ over $A\in H_2^{\mbox{\scriptsize{free}}}(X;\Z)$ where $c_A\in\Q$ and
the $t_A$ are
variables satisfying  $t_At_{B}=t_{A+B}$.

 After summing on $A$ and dualizing,
(\ref{modulimap2}) defines a map
\begin{equation}
GW_{g,n,X}:\, H^*(\overline{\M}_{g,n})\otimes H^*(X^n) \ \to \ RH_2(X).
\label{1.3}
\end{equation}
(This is well-defined for $2g+n\geq 3$, and we will always assume  that it
has been extended to
all $g,n$ by the stabilization procedure of [IP1] \S 3).

 Finally, we sum over $n$ and $g$ by setting  $\overline{\M} =
\bigcup_{g,n}\overline{\M}_{g,n}$,    letting ${\Bbb T}(X)$ denote the
total tensor algebra
${\Bbb T}(H^*(X))$ on the rational cohomology of $X$, and introducing a
variable  $\lambda$ to keep track of
the euler class.  The total Gromov-Witten invariant  are then the map
\bear
GW_X: H^*(\overline{\M})\otimes  {\Bbb T}(X)\ra RH_2(X)[\la].
\label{defGW}
\eear
which is linear in the second factor and defined by the Laurent expansion
\begin{equation}
GW_{X}\ =\ \sum_{g,n,A} \frac{1}{n!}\,GW_{g,n,A,X}\ t_A\ \la^{2g-2}.
\label{defn.7}
\end{equation}
The diagonal action of the symmetric group on
$\overline{\M}_{g,n}\ti X^n$ leaves $GW_X$ invariant up to sign, and if
$\kappa\in H^*(\overline{\M}_{g,n})$
then
$GW_X(\kappa,
\alpha)$ vanishes unless $\alpha$ is a tensor of length $n$.

We can recover the familiar geometric interpretation of these invariants
by  evaluating on cohomology classes.  Given $\kappa\in
H^*(\overline{\M};\Q)$ and  a
vector $\alpha=(\alpha_1,\dots, \alpha_n)$ of rational cohomology classes
in $X$ of length $|\alpha|=n$,
fix a generic $(J,\nu)$ and generic
geometric representatives $K$ and $A_i$ of the Poincar\'{e} duals of
$\kappa$ and of the $\alpha_i$. Then
$GW_{g,n,A,X}(\kappa,\alpha)$ counts (with orientation) the number of genus
$g$ $(J,\nu)$-holomorphic maps
$f:C\to X$ with $C\in K$ and  $f(x_i)\in A_i$ for each marked point $x_i$.  By the usual dimension counts,
this vanishes   unless
$$
\deg \kappa+\sum\deg \alpha_i- 2 |\alpha|=(\mbox{dim}\, X-6)(1-g)-2K_X\cdot A
$$
 where $K_X$ is the canonical class of $X$.

It is sometimes useful to also incorporate the `$\psi$ classes'.   There
are canonical oriented real 2-plane bundles ${\cal L}_i$ over
$\overline{\M}_{g,n}(X,A)$
whose fiber at each map $f$ is the  cotangent space to the domain curve
at the $i^{\mbox{\tiny th}}$ marked point.  Let $\psi_i$ be the euler class of
${\cal L}_i$, and for each vector $D=(d_1,\dots d_n)$ of non-negative
integers let $\psi_D=\psi_1^{d_1}\cup\dots \cup\psi_n^{d_n}$.  Replacing
the lefthand side of (\ref{modulimap2}) by the pushforward of the cap product
$\psi_D\cap \left[\overline{\M}_{g,n}(X,A)\right]$
and again dualizing gives invariants
\begin{equation}
GW_{g,n,D,X}:\, H^*(\overline{\M}_{g,n})\otimes H^*(X^n) \ \to \ RH_2(X)
\label{1.45}
\end{equation}
which agree with (\ref{1.3}) when $D$ is the zero vector.  These invariants
can be included in
(\ref{defGW}); the series (\ref{defn.7}) then has additional variables that keep
track of the vector $D$. To keep the notation manageable we will not do this
explicitly, although these
invariants will reappear at the end of section 4. (The similar classes
$\phi_i$ defined by the
 cotangent space of the stabilized domain are already part of (\ref{defGW})).

\bigskip

The  invariant (\ref{defGW}) counts
$(J, \nu)$-holomorphic maps from {\it connected}  domains.  It is often
more natural to work with domains
with more components, as  Taubes  does in [T] (see also [IP1]).  Let
$\wt{\M}_{\chi,n}$
be the space of all compact Riemann surfaces of euler characteristic $\chi$
with any number of components, of any possible genus, and with $n$ ordered
points distributed in all possible ways. More precisely, let ${\cal P}_n$ be
the set of all ordered partitions of the set $\{x_1,\dots x_n\}$. Each
$\pi=(\pi_1,\dots,\pi_{\ell})\in{\cal P}_n$ records which marked points go
on each component. Then
\best
\wt \M_{\chi,n}=\ma\bigsqcup_{\pi\in {\cal P}_n}\;
\ma\bigsqcup_{g_i}\  \overline{\M}_{g_1,\pi_1}\ti \dots \ti
 \overline{\M}_{g_l,\pi_l}
\eest
where $\overline{\M}_{g_i,\pi_i}$ is the Deligne-Mumford  space with $n_i$
marked points labeled by the set $\pi_i$ and where the second sum is over
all $g_i$ with
$\sum(2-2g_i)=\chi$.  Define the  ``Taubes-Witten'' invariant
$$
TW_X: H^*(\wt{\M})\otimes {\Bbb T}(X) \ \to \ RH_2(X)[\la]
$$
by
\begin{equation}
TW_X\ =\ e^{GW_X}.
\label{defnTW}
\end{equation}
As before, when this is expanded as a Laurent series
$$
TW_X(\kappa, \alpha)\ =\ \sum {1\over n!}\;TW_{X,A,\chi,n}(\kappa,\alpha)\
\la^{-\chi}\,t_A
$$
the coefficients count the number of curves (not necessarily connected)
with euler characteristic $\chi$ representing $A$ satisfying the
constraints $(\kappa,\alpha)$.

\medskip

\rem \ \ Specifically, if $\al=\al_1\otimes\dots \otimes \al_n$, then for
each partition $\pi=(\pi_1,\dots \pi_\ell)$
let $\al_{\pi_i}$ be the product of $\al_j$ for all $j\in \pi_i$. If
$\kappa=\kappa_1\otimes \dots \otimes \kappa_\ell$, then:
\best
TW_{X, n}(\kappa,\alpha)\ &=&\;\sum_{\pi\in{\cal P}_n}
{\ep(\pi) \over \ell! }{ n \choose n_1, \dots, n_\ell}
 GW_{X,n_1}(\kappa_{1},\al_{\pi_1})
\otimes  \dots \otimes GW_{X, n_\ell}(\kappa_{\ell},\alpha_{\pi_\ell})
\eest
where $\ep_\pi=\pm 1$ depending on the sign of the permutation
$(\pi_1,\dots, \pi_\ell)$ and the degrees of $\al$.
 Note that  $A$ and $\chi$ add when one
take disjoint unions, so the variables $t_A$ and $\la$ multiply.

\vskip.4in

%%%%%%%%%%%%%%%%%%%%%%%%%%%%%%%%%%%%%%%%%%%%%%%%%%%%%%%%%%%%m hj
%%%%%%%%%%%%%%%%%  Section 2 %%%%%%%%%%%%%%%%%%%%%%%%%%%%%%
%%%%%%%%%%%%%%%%%%%%%%%%%%%%%%%%%%%%%%%%%%%%%%%%%%%%%%%%%%%%%%

\setcounter{equation}{0}
\section{Relative Invariants}
\bigskip

The invariants of Section 1 can be extended to  invariants for $(X,\w)$
relative to a codimension
two symplectic submanifold
$V$  (a related invariant has been defined by Li-Ruan [LR]). This requires addressing three major issues:
``rim tori'',  higher order contact between the
holomorphic curves and $V$, and   ``chambers''.

\bigskip

Let ${\cal J}^V$ be the space of pairs $(J,\nu)$ where $J$ is an almost
complex structure on $X$
compatible with the symplectic form and for which $V$ is a holomorphic
submanifold, and $\nu$ is a
perturbation (as in [RT]) such that $\nu|_V$ takes values in the tangent
space to $V$. A
$(J,\nu)$-holomorphic map into $X$  is {\it $V$-regular} if the inverse image 
of $V$ is a finite set of
points.  For each $(g,n)$, the $V$-regular maps form an open subset
$\M_{g,n}^V(X)$ of the moduli space
 of maps, and this is a manifold for   generic $(J,\nu)\in{\cal J}^V$. We
will write this moduli space as a
disjoint union of strata and compactify each strata separately.

\medskip

The strata are  labeled by sequences $s=(s_1, s_2,\dots, s_\ell)$ of
integers with
$s_i\geq 1$.  Let ${\cal S}$ be the
set of all such sequences and define the degree, length, and order of
$s\in{\cal S}$ by
$$
\mbox{deg }s =\sum s_i \qquad \ell(s) = \ell, \qquad |s|=s_1 s_2\cdots s_\ell.
$$
Associated to each $s$ is the moduli space
$$
\M_{g,n,s}^V(X)\,\subset\ \M_{g,n+\ell(s)}^V(X)
$$
consisting of those maps $f$ such that $f^{-1}(V)$ is exactly the points
$x_i$, $n+1\leq i \leq n+\ell(s)$, each with multiplicity $s_i$. Then
$\M_{g,n,s}^V(X)$ defines a stratum of
$\M_{g,n}^V(X)$ via the projection that forgets the last
$\ell(s)$ points, and  $\M_{g,n}^V(X)$ is the disjoint union of these strata.

\medskip

 To compactify  $\M_{g,n,s}^V(X)$, we introduce the ``symplectic
compactification'' $\widehat{X}$ of $X\setminus V$.  This is a compact,
symplectic manifold  whose
boundary is the circle bundle $S$ of the normal bundle to $V$ in $X$
together with a projection
$\pi:\widehat{X}\to X$ which is the projection $S\to V$ on the boundary and
is a  symplectomorphism in the
interior. Let $V^*$ be $V$ with the discrete topology, and let $S^*$  be
$S$ topologized as the disjoint
union of its fiber circles.  Then
$\pi:S^*\to V^*$ and the inclusions $V^*\subset X$ and
$S^*\subset \hat{X}$ are continuous.  The long exact sequence of the pair
$(\hat{X},S^*)$ gives maps
\bear
H_2(\hat{X})  \longrightarrow H_2(\hat{X},S^*)
\overset{\rho}\longrightarrow H_1(S^*).
\label{2.1}
\eear

 Each $V$-regular map  lifts to  a map into $\hat{X}$; this determines an
element of $H_2(\hat{X},S^*)$
that projects under $\rho$ to  an element of the symmetric
product $\mbox{Sym}\,V$.  Similarly, each $f\in\M_{g,n,s}^V(X)$ lifts to
$\hat{X}$, and evaluation at the
last $\ell(s)$ marked points gives a map into  the space
$$
{\cal S}V\ =\ \bigcup_{s\in{\cal S}}\ V_s
$$
where each $V_s$ is a copy of $V^{\ell(s)}$.    We give $SV$ the
topology of that makes this a
disjoint union.  Restricting (\ref{2.1}) to $\mbox{Sym}\,V$ and pulling
back by the projection $
{\cal S}V \to \mbox{Sym}\,V$,  then gives a space  ${\cal
H}_X^V=H_2(\hat{X},S^*)
\ti_{\footnotesize{\mbox{Sym $V$}}} {\cal S}V$ with an evaluation map $\ep$
into ${\cal S}V$:
\begin{equation}\begin{array}{clcl}
H_2(\hat{X}) & \longrightarrow &{\cal H}_X^V & \\
& & \Big\downarrow\ep &  \\
& &  {\cal S}V. &
\end{array}
\label{HSVproj}
\end{equation}
 We topologize ${\cal H}_X^V$ so that $\ep$ is a covering map. With this
setup, lifting maps into
$\hat{X}$ gives a diagram
\begin{equation}\begin{array}{clccc}
\M_{g,n,s}^V(X) & \overset{\phi}\to & {\cal H}_X^V& \overset{\ep}\to &{\cal
S}V \\
\Big\downarrow & & \Big\downarrow & &\Big\downarrow\\
\M_{g,n}^V(X) &  \to &  H_2(\hat{X},S^*) & \overset{\rho}\to & \mbox{Sym}\, V.\\
\end{array}
\label{liftmap}
\end{equation}
Note that the intersection divisor with $V$ can be recovered by applying
$\rho$, and that the image of $\phi$ implicitly describes that class of the
lift in $H_2(\hat{X},S)$.
However, these  data cannot be separated because the covering
(\ref{HSVproj}) is not, in general, trivial.

\medskip

\begin{rem}
The covering (\ref{HSVproj}) can be partially trivialized.  Let ${\cal
R}_X$ be the kernel of the
projection $H_2(\hat{X})\to H_2(X)$.  Elements of ${\cal R}_X$ are called {\it rim tori} because they can
be represented by tori  of the form  $\alpha\ti \mu$ where $\alpha$ is  a
curve in $S$ and $\mu$ is a
fiber of $S\to V$.   Projecting 2-cycles in  $(\hat{X},S^*)$ into $X$
defines a map $\pi_*:{\cal H}_X^V\to
H_2(X)$.  Then (\ref{HSVproj}) can be rewritten as a covering map
\begin{equation}
\begin{array}{clc}
{\cal R}_X & \longrightarrow &{\cal H}_X^V \\
& & \Big\downarrow\\
& &   H_2(X)\ti {\cal S}V.
\end{array}
\label{Hfibration2}\nonumber
\end{equation}
\end{rem}

\medskip

\begin{lemma}
 For generic $(J,\nu)$, the image of $\M_{g,n,s}^V(X)$ under $\phi$  has
compact closure in ${\cal
H}_X^V$ for each $s\in {\cal S}$,  and these images define  elements of
$H_*({\cal H}_X^V;\Q)$.
\label{prop1}
\end{lemma}

The relative GW invariant of the pair $(X,V)$ is then obtained by
repeating the discussion leading from (\ref{modulimap2}) to (\ref{defGW}).
Writing ${\Bbb T}(H^*(X,V))$ as 
${\Bbb T}(X,V) $, the relative  invariant is thus  a map
$$
GW_X^V: H^*(\overline{\M})\otimes  {\Bbb T}(X,V) \ \longrightarrow\
 H_*({\cal H}_X^V)[\la]
$$
with an expansion
\bear
GW_X^V =\ \sum_{g,n} {1\over n!}\; GW_{g,n,X}^V\ \la^{2g-2}.
\label{defnRelInvt2}
\eear
 Formula (\ref{defnTW}) extends this to a relative Taubes-Witten invariant
\begin{equation}
TW_X^V:  H^*(\overline{\M})\otimes  {\Bbb T}(X,V)  \to    H_*({\cal
H}_X^V)[\la].
\label{defn5}
\end{equation}
When $V$ is the empty set ${\cal H}_X^V$ is $H_2(X)$ and the relative
invariant takes values
in $RH_2(X)$.  The relative and absolute invariants are then equal:
$GW_X^{\emptyset}=GW_X$.
\smallskip

The following theorem shows that these relative
invariants are well-defined, although they depend on a `chamber structure'
in ${\cal J}^V$. Let ${\cal W}$ be the subset  of ${\cal J}^V$ consisting
of all those
$(J,\nu)$ which admit a $(J,\nu)$-holomorphic map $C\subset V$ such that
the restriction of the linearization $D_C$ to the normal bundle
to $V$ has non-trivial cokernel.  This `wall' ${\cal W}$ has codimension
one and depends only on the 1-jet of
$(J,\nu)$ along $V$.  The wall separates ${\cal
J}^V\setminus{\cal W}$ into components called  the chambers of ${\cal J}^V$.

\begin{theorem}
\begin{tabular}[t]{l}
 $GW^V$ and $TW^V$ are constant for generic $(J,\nu)$ in each chamber of
${\cal J}^V$.
\end{tabular}
\label{thm2.3}
\end{theorem}

\medskip

 The geometric interpretation of the relative invariant
(\ref{defnRelInvt2}) is similar to that of
(\ref{defn.7}).
Given $\kappa\in H^*(\overline{\M})$, a vector
$\alpha=(\alpha_1,\dots, \alpha_n)$ of classes in $H^*(X,V)$, and a
$\gamma\in H^*({\cal H}_X^V)$, fix generic
geometric representatives $K\subset \overline{\M}_{g,n}$,
$\Gamma\subset {\cal H}_X^V$ and $A_i\subset X$ of their Poincar\'{e}
duals. Then the evaluation pairing
$\gamma \cdot GW_{g,X}^V(\kappa,\alpha)$ counts   the oriented number of
genus $g$ $(J,\nu)$-holomorphic maps $f:C \to X$ with $C\in K$,
$\phi(f)\in \Gamma$, and $f(x_i)\in A_i$ for each marked point $x_i$.

Note that the condition $\phi(f)\in \Gamma$ constrains both the homology
class $A$ of the map and the boundary values of the curve.  In the special
case when there are no rim
tori,  these homology and the boundary value constraints can be fully separated.

\medskip

\begin{ex} When there are no rim tori, ${\cal H}_X^V$ is the subset of
$H_2(X)  \ti \bigcup_s V^{\ell(s)}$
consisting of pairs $(A,s)$ with $\deg s=A\cdot V$.  The homology of
${\cal H}_X^V$ is the corresponding
subalgebra of $RH_2(X)\otimes {\Bbb {CT}}$ where ${\Bbb {CT}}$ is the
``contact tensor algebra'' of $V$:
\best
{\Bbb {CT}}(V)\,=\, {\Bbb {T}}(H_*(V)\ti {\Bbb N}).
\eest
The relative invariants are then maps
\best
 H^*(\overline{\M})\otimes {\Bbb T}(X,V) \ra  {\Bbb {CT}}(V)\otimes RH_2(X)[\la]
\eest
and have Laurent expansions like (\ref{defn.7}) with coefficients in ${\Bbb
{CT}}(V)$.  In fact, a basis for
the contact algebra is given by  elements of the form
$C_{s,\gamma}=C_{s_1, \gamma_1}\otimes\dots \otimes  C_{s_\ell, \gamma_\ell}$,
where $\gamma_i$ are elements of $H_*(V;\Q)$ and $s_i\ge 1$ are integers.
Let $\{C^*_{s,\gamma}\}$ denote the
dual basis.  With $\kappa$ and $\alpha$ as above, we can expand
\bear
 \sum_{s,\gamma}  GW_{g,A,X}^V(\kappa,\alpha; C_{s,\gamma})\ \
C^*_{s,\gamma}\ t_A\ \la^{2g-2}
\label{2.last}
\eear
where the coefficients count the oriented number of genus $g$
$(J,\nu)$-holomorphic
maps $f:C \to X$ with $C\in K$, $f(x_i)\in A_i$, and having a contact of
order $s_j$ with $V$ along fixed representatives $\Gamma_j$ of the
Poincar\'{e} duals of the $\gamma_j$.
\label{ex2.3}
\end{ex}

\vskip.4in

%%%%%%%%%%%%%%%%%%%%%%%%%%%%%%%%%%%%%%%%%%%%%%%%%%%%%%%%%%%%
%%%%%%%%%%%%%%%%%  Section 3 %%%%%%%%%%%%%%%%%%%%%%%%%%%%%%
%%%%%%%%%%%%%%%%%%%%%%%%%%%%%%%%%%%%%%%%%%%%%%%%%%%%%%%%%%%%%%

\setcounter{equation}{0}
\section{The Convolution}
\bigskip

 Having defined the relative invariants, we now turn to the problem of
describing how they behave under
the symplectic  sum operation. The goal is to give a formula
expressing the absolute invariants of a
symplectic sum $Z=X\#_{V=\overline{V}} Y$ in terms of the relative
invariants of $X$ and
$Y$.  This turns out to be most
natural when one works with the $TW$ invariants, and when one works with
relative invariants throughout.

\bigskip

Let $X$ be a symplectic manifold with two disjoint symplectic  submanifolds
$U$ and $V$ with real
codimension two.  Suppose that $V$ is symplectically identified with a
submanifold of
similar triple  $(Y,\overline{V},W)$ and that the normal bundles of
$V\subset X$ and
$\overline{V}\subset Y$ have opposite chern 
classes.  Let $(Z,U,W)$ be the resulting symplectic  sum.  We will choose
$(J,\nu)$ on the $X$ and $Y$
sides so that the chambers along $V$ and $\overline{V}$  match when the sum
is formed.  In  this section we
will define a ``convolution'' operation  and establish a formula of the  form
\bear
TW_X^{U\cup V}\, *\, TW_Y^{ \overline{V}\cup W}\, = \, TW_Z^{U\cup W}
\label{3.1}
\eear
under the assumption that all curves contributing to the invariants are
$V$-regular (this condition will be
eliminated in the next section).

\medskip

The convolution operation is assembled from several pieces.   The first
ingredient is a gluing map defined
on Deligne-Mumford spaces.  Given stable curves $C_1$ and
$C_2$ (not necessarily connected) with euler
characteristics $\chi_i$ and $n_i+\ell$ marked points, we can construct a
new curve by identifying the last
$\ell$ marked points of
$C_1$ with the last $\ell$ marked points of $C_2$, and then
forgetting the marking of these new
double points.  This  defines an attaching  map
\best
\xi_\ell:\wt{\M}_{\chi_1,n_1+\ell}\ti \wt{\M}_{\chi_2,n_2+\ell}\ \longrightarrow
\ \wt{\M}_{\chi_1+\chi_2-2\ell,n_1+n_2}
\eest
whose image is a subvariety of complex codimension $\ell$.  Taking the union
over all
$\chi_1, \chi_2, n_1$ and $n_2$ gives a  map
$\xi_\ell:\wt{\M}\ti\wt{\M}\to\wt{\M}$ for each $\ell$.

The second ingredient describes how maps are glued along $V$.  Consider the
evaluation map
$$
\ev: \M^V(X) \ti \M^{\overline{V}}(Y)  \ \overset{\phi}\longrightarrow \
{\cal H}_X^V
\ti {\cal H}_Y^{\overline{V}}
\ \overset{\ep}\longrightarrow\  {\cal S}V\ti {\cal S}V.
$$
defined by the top row of (\ref{liftmap})  on each factor;
this evaluation map
 records the intersection points with $V$, together with their
multiplicities.  The diagonal of
${\cal S}V\ti {\cal S}V$ is the union of components
$$
{\bold \Delta}_s\ \subset \ {\cal S}V\ti {\cal S}V.
$$
labeled by sequences $s$.   With this notation, we can state our first main
result.  It is a `gluing
theorem' for families of pseudo-holomorphic maps which is proved by PDE
methods [IP2].

\begin{theorem}  For each sequence $s$, there is an
${|s|\over \ell(s)!}$-fold cover of $ev^{-1}({\bold \Delta}_s)$,
denoted \newline $\M^V(X) \ti_{ev_s} \M^{\overline{V}}(Y)$, and a commutative
diagram
$$
\begin{array}{ccc}
\M^V(X) \ti_{ev_s} \M^{\overline{V}}(Y) & \longrightarrow &
\M(Z)\\
\Big\downarrow\vcenter{\rlap{$st$}} &  & \Big\downarrow
\vcenter{\rlap{$st$}}  \\
\wt{\M}\ti \wt{\M} &
\overset{\xi_{\ell(s)}}\longrightarrow &  \wt{\M}
\end{array}
$$
where the top arrow is a diffeomorphism.
\label{gluing thm}
\end{theorem}

\medskip

The next step is to pass to homology.  We can glue a curve in $X$ to a
curve in $Y$ provided the two
curves meet $V$ at the same points with the same
multiplicity; this defines a map
\bear
g:{\cal H}_X^{U,V}\ti_{\ep} {\cal H}_Y^{\overline{V},W}\ra {\cal H}_Z^{U,W}.
\label{gmap}
\eear
\non Combining this with the map $\xi_\ell$ above gives the convolution
operator that describes how homology
classes of maps combine in the gluing operation for the symplectic sum.

\begin{defn}\ \ The convolution operator
$$
*: \ H_*(\wt{\M}\ti {\cal H}_X^{U,V};\Q[\la])
\otimes H_*(\wt{\M}\ti {\cal H}_Y^{\overline{V},W};\Q[\la])
\ \longrightarrow\  H_*(\wt{\M}\ti {\cal H}_Z^{U,W};\Q[\la])
$$
is given  by
\begin{equation}
(\kappa\otimes h)\,*\,(\kappa'\otimes h') \ =\
\sum_{s} {|s|\over \ell(s)!} \
\left(\xi_{\ell(s)}\right)_*(\kappa\otimes \kappa')  \
g_*\left[ \left. h\ti h'\right|_{ \ep^{-1}\left(\bold
{\Delta}_{s}\right)}\right]\ \la^{2\ell(s)}.
\label{wteddiagonal1}
\end{equation}
\end{defn}

\medskip

 To conform with this definition, we will regard the
$TW$
invariant of
$X$ as a map
$$
TW_X^{U,V}:    {\Bbb T}(X,V)  \to   H_*(\wt{\M}\ti {\cal H}_Z^{U,V};\Q[\la])
$$
 by dualizing the first term in (\ref{defn5}). Finally, we say that a
constraint $\alpha\in H^*(Z)$ {\it separates} as $\alpha_X+\alpha_Y$ if its
Poincar\'{e} dual
$PD(\alpha)$ is the image of $(PD(\alpha_X),PD(\alpha_Y))$ under the
Mayer-Vietoris map $H_*(X\setminus
V)\oplus H_*(Y\setminus \overline{V})\to H_*(Z)$. These generate a subalgebra
${\Bbb T}_s(Z)$ of ${\Bbb T}(Z)$,  namely the span
of tensors of the form
\begin{equation}
\label{separates}
\alpha= \alpha_X+\alpha_Y
\end{equation}
with  $\alpha_X\in
{\Bbb T}(X,V)$ and
$\alpha_Y\in {\Bbb T}(Y,\overline{V})$.  Then, when $\alpha$ has this form
and assuming that all
curves contributing to the
invariants are $V$-regular,  (\ref{3.1}) holds in the sense that
\bear
TW_Z^{U, W}(\alpha)\, = \,TW_X^{U,V}(\alpha_X)\, *\, TW_Y^{
\overline{V}, W}(\alpha_Y).
\label{3.5}
\eear

\bigskip

\begin{ex} {\bf (The S-matrix)}\ \  Starting from the normal bundle $N_V$
of $V$ in $X$, we can form the
$\P^1$ bundle ${\Bbb F}=\P(N_V\oplus\cx)$ over $V$.  In  ${\Bbb F}$, the
zero section $V$   and the infinity section $\overline{V}$ are disjoint
symplectic submanifolds.  The
symplectic  sum of $(X,U,V)$ and $({\Bbb F}, \overline{V},V)$  is a
symplectic deformation of $(X,U,V)$,
so has the same  $TW$ invariant. The convolution then defines a operation
$$
 \ H_*(\wt{\M}\ti {\cal H}_X^{U,V};\Q[\la])
\otimes H_*(\wt{\M}\ti {\cal H}_{{\Bbb F}}^{\overline{V},V};\Q[\la])
\ \longrightarrow\  H_*(\wt{\M}\ti {\cal H}_X^{U,V};\Q[\la]).
$$
Thus for each choice of constraints $\alpha\in {\Bbb T}({\Bbb F},V\cup
\overline{V})$, the
$TW$ invariant of ${\Bbb F}$ relative to its zero and infinity section
defines an ``S-matrix''
\bear
\label{3.Endo}
TW_{{\Bbb F}}^{\overline{V},V}(\alpha)\,\in\, \mbox{End}\
\left(H_*(\wt{\M}\ti {\cal H}_X^{U,V};\Q[\la])  \right)
\eear
which describes how families of curves on $X$
are modified --- ``scattered''--- as they pass through a neck modeled on
$({\Bbb F}, \overline{V},V)$  containing the constraints $\alpha$.
\label{ex3.1}
\end{ex}

\bigskip

The identity  in the endomorphism algebra (\ref{3.Endo}) is always realized
as  the convolution
by the element
\best
{\Bbb I}\in H_*(\wt{\M}\ti {\cal H}_{{\Bbb F}}^{\overline{V},V};\Q[\la])
\eest
corresponding to that part of $TW(1)$ coming from the maps whose domain is
unstable and whose
image is a multiple cover of a fiber of ${\Bbb F}$.  As in physics, it is
convenient to write the $S$-matrix
(\ref{3.Endo}) as ${\Bbb I}+ R^{\overline{V},V}$; the inverse of the
$S$-matrix is then given by the series
\bear
(TW_{\Bbb F}^{\overline{V},V})^{-1}={\Bbb I}\ -\  R^{\overline{V},V} \ +
\  R^{\overline{V},V}*  R^{\overline{V},V}\ -\
R^{\overline{V},V}* R^{\overline{V},V}*R^{\overline{V},V}\ +\ \dots.
\label{3.inverseSmatrix}
\eear

\bigskip

\begin{ex}
When $V=\P^1$,   ${\Bbb F}\to V$ is one of the rational ruled surfaces
with its standard symplectic and holomorphic structure.  If we wish to
count all pseudo-holomorphic maps,
without constraints on the genus
or the induced complex structure, the relevant S-matrix is  the relative
$TW$ invariant  with
$(\kappa,\alpha)= (1,1)\in H^0(\overline{\M})\ti H^0(X)$.  This  case
works out
neatly:  a dimension count shows that  $TW_{{\Bbb
F}}^{\overline{V},V}(1,1)={\Bbb I}$.
\label{ex3.2}
\end{ex}

\vskip.4in

%%%%%%%%%%%%%%%%%%%%%%%%%%%%%%%%%%%%%%%%%%%%%%%%%%%%%%%%%%%%
%%%%%%%%%%%%%%%%%  Section 4 %%%%%%%%%%%%%%%%%%%%%%%%%%%%%%
%%%%%%%%%%%%%%%%%%%%%%%%%%%%%%%%%%%%%%%%%%%%%%%%%%%%%%%%%%%%%%

\setcounter{equation}{0}
\section{The Sum Formula}
\bigskip

The full  formula for the invariants of a symplectic sum is obtained from
(\ref{3.5}) by a `stretching the
neck' argument.   The  sum $X\#_VY$ can be viewed as assembled from
$n+1$ symplectic sums, where the
middle $n$ pieces are copies of the ruled space ${\Bbb F}$ associated to
$V$. We again pick  $(J,\nu)$ on the pieces of this sum so that the
chambers match along each
cut.  The  pigeon-hole principle implies that, for each $A\in H_2(Z)$, we
can take $n$ sufficiently large so
that  all pseudo-holomorphic  $A$-curves in $Z$ limit to $V$-regular curves
in the $n+1$-fold
 sum.  Thus the
curves in $Z$ decompose into curves in $X$ joined to curves in $Y$ by a
chain of curves in the
intermediate copies of ${\Bbb F}$. Applying formula (\ref{3.5})  along each
cut yields the general sum formula.

\medskip

\begin{theorem}  Let $(Z,U,W)$ be the symplectic  sum of
$(X,U,V)$ and  $(Y,\overline{V},W)$ along $V=\ov V$.  Suppose that 
$\alpha\in{\Bbb T}(Z)$
separates as in
(\ref{separates}).  Then the relative TW invariant of $Z$ is given
in terms of the  invariants of $(X,U,V)$ and $(Y,\overline{V},W)$ by
$$
TW_Z^{U,W}(\al)\ =\   TW_X^{U,V}(\al_X) \, * \,
(TW_{\Bbb F}^{\overline{V},V})^{-1}\, * \, TW_Y^{\overline{V},W}(\al_Y)
$$
where  the middle term is the  inverse of the $S$-matrix
(\ref{3.inverseSmatrix}).
\label{mainthm}
\end{theorem}

One interesting special case is when $U$ and $W$ are empty and we put no
constraints on the
complex structure of the curves.  Theorem \ref{mainthm} then gives a
formula expressing the absolute
invariant of $Z$ in terms of the relative invariants of $X$ and $Y$.

In this case the algebra of the formula simplifies:  we drop all references 
to $\wt{\M}$ and, because $U$ and
$W$ are  empty, $H_*({\cal H}_Z^{U,W})$ is simply
$RH_2(Z)$.  The convolution (\ref{wteddiagonal1}) can then be thought of as
a bilinear pairing

$$
\langle\ ,\ \rangle:\ H_*({\cal H}_X^V;\Q[\la])
\otimes H_*({\cal H}_Y^{\overline{V}};\Q[\la])
\longrightarrow  RH_2(Z)[\la].
$$

\medskip

\non In this case $g:{\cal H}_X^{V}\ti_{\ep} {\cal H}_Y^{\overline{V}}\ra
RH_2(Z)$, so for each sequence $s$, 
$\ep^{-1}\left({\bold \Delta}_s\right)$ is the union of components ${\bold
\Delta}_{A,s}=\ep^{-1}\left({\bold \Delta}_s\right) \cap g^{-1}(A)$,  and
the formula for the bilinear
pairing becomes
$$
\langle  h\,, \, h'  \rangle
\ =\  \sum_{s} {|s|\over \ell(s)!}  \sum_{A\in H_2(Z)}\ g_*[{\bold
\Delta}_{A,s}\cap (h\ti h')]\ t_A \
\la^{2\ell(s)}.
$$
The $S$-matrix is  the operator  obtained by letting
$TW_{{\Bbb F}}^{\overline{V},V}(1,1)\,\in\, H_*({\cal H}_X^{U,V};\Q[\la])
$ act  by  convolution; let $S_{\overline{V},V}$ denote its inverse.

\medskip

\begin{cor}  Let $Z$ be the symplectic  sum of $(X,V)$ and
$(Y,\overline{V})$ and suppose that $\alpha\in{\Bbb T}(Z)$ separates as in
(\ref{separates}). Then
$$
TW_Z(\al)\ =\ \left\langle\, TW_X^V(\al_X),\ S_{\overline{V},V}\,
\cdot\,TW_Y^{\overline{V}}(\al_Y)
\,\right\rangle.
$$
\label{cor4.2}
\end{cor}

Another simplification occurs when there are no rim tori in $X$ and $Y$,
and therefore in $Z$.
Then  the
relative invariants have an expansion of the form (\ref{2.last}).  For
simplicity we take $U$ and $W$ both
empty.  Then ${\cal H}_Z^{U,W}$ is simply  $H_2(Z)$ and the map $g$ of
(\ref{gmap}) is the
identification ${\cal S}V\ti_\rho {\cal S}V$ with ${\cal S}V$ together with the addition
$H_2(X)\ti H_2(Y)\to H_2(Z)$.  The
key part of the convolution (\ref{wteddiagonal1}) is then given by the cap
product with the Poincar\'{e}
dual of the diagonal:
$$
g_*\left[ \left.h\ti h'\right|_{\bold {\Delta}_{s}}\right]\ =\
\sum_{B\in H_2(X)}\sum_{ C\in H_2(Y)} {\mbox{PD}}\left(\bold
{\Delta}_{s}\right)\cap (h\ti h')\ t_B t_C.
$$
We can then `split the diagonal' by fixing a basis $\{h^{p}\}$ of
$H^*({\cal S}V)$
and writing
$$
{\mbox{PD}}\left(\bold {\Delta}_{s}\right)\ =\
\sum_{p,q} Q^V_{p,q}\ h^p\ti h^q= \
\sum_{p} \ h^p\ti h_p
$$
where $Q^V_{p,q}$ is the intersection matrix of the Poincar\'{e}
dual basis, and $\{h_p\}$ is the dual basis with respect to $Q^V$.  In fact,
fix a basis $\{\gamma^{i}\}$ of $H_*(V)$, and let
$\{\gamma_{i}\}$ be
the dual basis of $H_*(V)$ with respect to the intersection form.
 Then the $\{h^{p}\}$ can be taken to be the
basis $\{C^*_{s,\gamma^I}=C^*_{s,\gamma^1}\otimes\dots \otimes
C^*_{s_\ell,\gamma^\ell}\}$ of  Example 2.4.
The convolution then has the more explicit form

$$
(\kappa\otimes h)\,*\,(\kappa'\otimes h') \ =\
\sum_{B,C,s} {|s|\over \ell(s)!} \
\left(\xi_{\ell(s)}\right)_*(\kappa\otimes \kappa')  \
 \sum_{I} C^*_{s,\gamma^I}(h)\,  C^*_{s,\gamma_J}(h') \
t_{B+C}\la^{2\ell(s)}.
$$

Finally, we can also include the $\psi$-classes (\ref{1.45}) to obtain a
sum formula for the
numbers
$$
TW_{\chi,A,D,X}(\kappa,\alpha)
$$
in terms of the relative invariants (also numbers)
$TW^V_{\chi,A,D,X}(\kappa,\alpha ;
C^*_{s,\gamma})$.    For this, we will take $\kappa=1$ and  pair  each
$\psi$ class with an $\alpha$.
This case is usually written in terms of the classes $\tau_d(\alpha)=
\psi^d_i\cup ev_i^*(\alpha)$, so we adopt the notation
$$
TW^V_{\chi,A,X}\left(\tau_D(\alpha) ;
 s,\gamma\right)\ :=\ TW^V_{\chi,A,D,X}(1,\alpha ; C^*_{s,\gamma}).
$$
Using dual bases $\{\gamma^I,\gamma_I\}$ as above, we let
$TW_{\chi, A}(s_1,\gamma^{I};s_2,\gamma_{J})$ denote the relative invariant
of ${\Bbb F}$
satisfying the contact constraints $ C^*_{s_1,\gamma^{I}}$ along
$\overline{V}$ and $C^*_{s_2,\gamma_{J}}$ along $V$.    Then the inverse
$S$-matrix is the alternating
series
\bear
S^{\overline{V},V}_{\chi, A}(s_1,\gamma_{I}; s_2,\gamma^{J})
 & = & {\Bbb I}(s_1,\gamma_{I}; s_2,\gamma^{J})\ -\
          TW^{\overline{V},V}_{\chi, A}( s_1,\gamma_{I}; s_2,\gamma^{J})\\
& + &\hskip-.2in
  \sum_{A_1+A_2=A\atop\chi_1+\chi_2-2\ell(s_3)=\chi} \ma\sum_{s_3,
I_3}{|s_3|\over \ell(s_3)!}\
TW^{\overline{V},V}_{\chi_1, A_1}(s_1,\gamma_I; s_3,\gamma^{I_3})\,
TW^{\overline{V},V}_{\chi_2, A_2}(s_3,\gamma_{I_3}; s_2,\gamma^{J})-\dots
\nonumber
\label{4.longprop}
\eear

\begin{theorem} Suppose there are no rim tori in $(X,V)$ and
$(Y,\overline{V})$ and that the constraint $\al$ separates as in
(\ref{separates}). Then the number
$$
TW_{\chi,A,Z}(\tau_D(\alpha))
$$
is given by
\best
\sum  \  \ \frac{|s_1|\, |s_2|}{\ell(s_1)!\, \ell(s_2)!}\ \
TW^V_{\chi_1,A_1,X}(\tau_{D}(\alpha_X);s_1,\gamma^I)\
S^{\overline{V},V}_{\chi_0,A_0}(s_1,\gamma_{I};s_2,\gamma^{J})
\ TW^{\overline{V}}_{\chi_2,A_2,Y}(s_2,\gamma_J;\phantom{i^I}
\tau_{D}(\alpha_Y))
\label{easymainformula}
\eest
where the sum is over all (a) decompositions $A=A_1+A_0+A_2$ and
$\chi=\chi_1+\chi_0+\chi_2-2\ell(s_1)-2\ell(s_2)$, (b) 
pairs of sequences $s_1,s_2\in{\cal S}$ of degree $\deg s_i=
A_i\cdot V$, and (c) pairs of multiindices $I,J$  enumerating  bases
$\{\gamma^I,\gamma_I\}$ of
$H^*(V^{\ell(s_i)})$ dual by the intersection form.
\label{2easymainformula}
\end{theorem}

In the case where $TW_{V,\overline{V}}(1,1)={\Bbb I}$ only the leading term
of (\ref{4.longprop}) is present, and
this sum formula   simplifies to
$$
\sum_{A=A_1+A_2\atop
\chi_1+\chi_2-2\ell(s)=\chi}\sum_{s,I}
 \frac{|s|}{l(s)!}\
TW^V_{\chi_1,A_1,X}\left(\tau_{D_X}(\alpha_X); s,\gamma^I\right)
\ TW^V_{\chi_2,A_2,Y}\left(s,\gamma_I ;\tau_{D_Y}(\alpha_Y)\right).
$$

\vskip.4in

%%%%%%%%%%%%%%%%%%%%%%%%%%%%%%%%%%%%%%%%%%%%%%%%%%%%%%%%%%%%m hj
%%%%%%%%%%%%%%%%%  Section 5 %%%%%%%%%%%%%%%%%%%%%%%%%%%%%%
%%%%%%%%%%%%%%%%%%%%%%%%%%%%%%%%%%%%%%%%%%%%%%%%%%%%%%%%%%%%%%

\setcounter{equation}{0}
\section{Applications}
\bigskip

This section presents two consequences of Theorem \ref{mainthm}:  the
Caporaso-Harris recursion formula for
the number of curves in $\P^2$, and the  formula (related to the
``Yau-Zaslow conjecture'') for the  number
of curves in the rational elliptic surface.  While these formulas
are both known, the proofs outlined here are considerably easier
and more transparent than the existing ones.

\bigskip

We begin by examining the relative invariants for the pair $(\P^2,L)$ (the
manifold $\hat{X}$ of section 2
is then the  4-ball, and in that sense we are considering the  relative
invariant
of $(B^4,S^3)$).  This pair can
be written as a symplectic  sum
$$
(\P^2,L)\ma\#_{L=E} ({\Bbb F}_1,E,L)\ = \ (\P^2,L)
$$
where $({\Bbb F}_1,E,L)$ is the ruled surface with euler class one with
its zero section $L$ and its infinity section $E$.  We will obtain a
recursive formula for the relative invariants of $(\P^2,L)$ by moving one point
constraint $p$ to the ${\Bbb F}$ side, and applying the sum formula.

 Following [CH], let
$N^{d,\delta}(\al,\b)$ be the number of degree $d$ curves in $\P^2$  with
$\delta$ double points with a
contact with $L$ of order $k$ at $\al_k$ fixed  points, and at $\b_k$
moving points, and passing through the
appropriate number of fixed generic points off $L$.   Similarly, let
$N^{a,b,\delta}_{\Bbb
F}(\al',\b';\; p\; ; \al,\b)$ be the number of curves in ${\Bbb F}$
representing $aL+bF$
that have contact described by
$(\al', \b')$ along $E$, $(\al, \b)$ along $L$, and pass through  $p\in {\Bbb
F}$.  By elaborating on Example 2.4, one can relate these numbers $N^{d,\delta}$
 to the relative $TW$ invariants.   Theorem \ref{mainthm} then becomes
\bear
 N^{d,\delta}(\al,\b)=\ma\sum   |\al'|\,|\b'|\  \
 N^{d',\delta'}(\al',\b') \ \,
N^{d-d',d,\delta''}_{\Bbb F}(\b',\al';\; p\; ; \al,\b)
\label{5.CH1}
\eear
(the $S$-matrix is the identity by Example  3.4).  Furthermore, a
dimension count shows that   there  are exactly two types of curves that
contribute to the
last term in (\ref{5.CH1}):  either (i)  several rational covers of the
fiber, one of them passing
through the point $p$, or (ii) several rational covers of the fiber and a
rational curve in the class $L+a F$ passing through $p$ and having all
contact points with $E$ and $L$ fixed.  In each situation there is a unique
such curve. In the case (i) nothing changes, except that a contact of some
order $k$
at a moving point becomes a contact of order $k$ at the fixed point $p_0$
where the fiber through $p$ hits $E$; in case (ii) the degree of the curve
drops by 1, and the curve looses some of the contact conditions at fixed
points and gains some at moving points, so $\al\ge \al'$ and  $\b'\ge \b$.
Thus the product in (\ref{5.CH1}) becomes a
sum of two terms.  The resulting formula has the form
\best
\ N^{d,\delta}(\al,\b)= \ma \sum k
N^{d,\delta'}(\al-\ep_k,\b+\ep_k) +
\ma \sum |\b'-\b|\cdot {\al\choose \al'} {\b'\choose b}
\cdot N^{d-1,\delta'}(\al',\b')
\eest
and is exactly the formula of Caporaso-Harris ([CH]  Theorem 1.1).

\smallskip

The same method, applied to the sum
${\Bbb F}_n\#{\Bbb F}_n={\Bbb F}_n$ yields the recursion formula of Vakil [V].

\bigskip

We next consider the  invariants of the rational elliptic surface
$E$.  Let $f$ and $s$ denote, respectively, the homology class of the fiber
and the section of the
elliptic fibration $E\to\P^1$.  We will consider the invariants
$GW_{g,s+df}$, which count the number of
connected genus $g$ curves in the class $s+df$ through $g$ generic points.
These define  power
series
$$
F_g(t)=\ma\sum_{d\ge 0} GW_{g,s+df}(p^g)\  t^{d}\,t_s
$$
where $t=t_f$. Recently,  G\"{o}ttsche  and Pandharipande [GP] obtained
recursive relations for the
coefficients of $F_0$,  and Bryan-Leung [BL]  proved the  closed-form expression
\bear
F_0(t)=t_s\, \left(\ma\prod_{d}{1 \over 1-t^d}\right)^{12}.
\label{BL}
\eear
(this formula was very strongly motivated by the work of Yau-Zaslow [YZ]).   We will prove this by relating
it to  the similar series of elliptic ($g=1$)  invariants
\best
H(t)&=&\ma\sum_{d\ge 0}\ GW_{1,s+df}(\tau_1(f^*)) \  t^{d}\,t_s
\eest
where $f^*\in H^2(E)$ is the Poincar\'{e} dual of the fiber class. The
``topological recursion relation''
for $g=1$ (TRR)  easily implies that
\bear
H(t)\ =\ \frac{1}{12} \left(t\,F_0' -F_0\right) +F_0 G
\label{rel1}
\eear
where $G=\sum\sigma(d)t^d$ is the generating function for
$\sigma(d)=\sum_{a|d}a$.   On the other hand, we
can write $E=E\#_{F}T^2\ti S^2$, and compute $H$ by moving the
constraint onto the $T^2\ti S^2$
side and using Theorem \ref{mainthm}. For this situation, the only possible
splitting is a genus 0
curve on $E$ side representing $s+d_1F$ and a genus 1 curve on the
$T^2\ti S^2$ side. But the genus 1 invariant of $T^2\ti S^2$
is computable from TRR, and  comes out to be $2(G(t)-1/24)$. Therefore,
$$
H(t)\ =\ 2 F_0\l(G-{1\over 24}\r).
$$
Equating this with (\ref{rel1}), we see that $F_0$ satisfies the ODE
\best
t\, F_0'=12\,G F_0
\eest
with $F_0(0)=GW_{0,s}t_s = t_s$.  Integrating and rearranging  gives (\ref{BL}).

\vskip.4in

\small

%%%%%%%%%%%%%%%%%%%%%%%%%%%%%%%%%%%%%%%%%%%%%%%%%%%%%%%%%%%%%%%%%%%%%%%%%%%%%%%

\end{document}